\documentclass[12pt,reqno,intlimits]{article}

\usepackage{amssymb,amsthm,amsmath,epsfig,graphicx}

\theoremstyle{remark}

\begin{document}

\begin{center}
{\Large {\bf INTEGRAL OPERATORS IN GRAND MORREY SPACES}}
%\pagenumbering{arabic} \setcounter{section}{3}
\end{center}
\vskip+1cm

\begin{center}
{\bf Alexander Meskhi}
%\pagenumbering{arabic} \setcounter{section}{3}
\end{center}
\vskip+1cm

{\bf Abstract.} We introduce grand Morrey spaces and establish the
boundedness of Hardy--Littlewood maximal, Calder\'on--Zygmund and
potential operators  in these spaces. In our case the operators
and grand Morrey spaces are defined on quasi-metric measure spaces
with doubling measure. The results are new even for Euclidean spaces.

\vskip+0.5cm

\vskip+0.4cm


\vskip+1cm
\section*{Introduction}

In  the paper we introduce the grand Morrey spaces
$L^{p),\theta,\lambda}$ and derive the boundedness of a class of
integral operators (Hardy--Littlewood maximal functions,
Calder\'on--Zygmund singular integrals and potentials) in these
spaces. We study the boundedness problem in the frame of
quasi-metric measure spaces with doubling measure but the results
are new even for Euclidean spaces.

The classical grand Lebesgue spaces $L^{p)}$ were introduced in
the paper by T. Iwaniec and C. Sbordone \cite{IwSb} when they
studied the problem of the integrability of the Jacobian $J(f,x)$
of the order preserving mapping $f = (f_1, \cdots, f_n): \Omega
\to {\Bbb{R}}^n$ under minimal hypothesis, where $\Omega$ is a
bounded domain in ${\Bbb{R}}^n$ and $n\geq 2$.

Later the generalized grand lebesgue spaces $L^{p), \theta}$
appeared in the paper  by L.  Greco, T.  Iwaniec  and C.  Sbordone
\cite{GrIwSb}, where the existence and uniqueness of the nonhomogeneous $n$- harmonic equation $div \; A(x, \nabla u)= \mu$ were established.

Structural properties of these spaces were investigated in the
papers \cite{Fi}, \cite{FiKa}, \cite{CaFi} etc.

A. Fiorenza, B. Gupta and P. Jain \cite{FiGuJa} proved the
boundedness of the Hardy--Littlewood maximal operator defined on
an interval in weighted $L^{p)}$ space, while the same problem for
the Hilbert transform and other singular integrals were studied in
the papers \cite{KoMe},  \cite{Ko}.

The Morrey spaces $L^{p, \lambda}$, which were introduced by C. Morrey in 1938 (see
\cite{Mo}) in order to study regularity questions which appear in
the calculus of variations, describe local regularity more
precisely than Lebesgue spaces and widely use not only harmonic
analysis but also partial differential equations (c.f. \cite{Gia},
\cite{GiTr}).

For essential  properties of $L^{p, \lambda}$ spaces and the
boundedness of maximal, fractional and singular operators in these
spaces we refer to the papers \cite{Ad1}, \cite{Pe}, \cite{ChiFr},
\cite{DiRa}, \cite{SaTa},  \cite{Gia}, etc.

Finally we mention that necessary and sufficient conditions for
the boundedness of maximal operators and Riesz potentials in the
local Morrey--type spaces were derived in \cite{BG}, \cite{BGG}.

%By the symbol $p'$ we denote the conjugate number of $p$, i.e. $p'
%:= \frac{p}{p-1}$, $1<p<\infty$.

\section{Preliminaries}
%\section{Maximal Operator in Grand Morrey Spaces}
\setcounter{equation}{0}

Let $X:=(X, \rho, \mu)$ be a topological space with a complete
measure $\mu$ such that the space of compactly supported
continuous functions is dense in $L^1(X,\mu)$ and there exists a
non-negative real-valued function  (quasi-metric) $d$ on $X\times
X$ satisfying the conditions:

\vskip+0.1cm

\noindent(i) $\rho(x,y)=0$ if and only if  $x=y$;

%\noindent(ii) $d(x,y)> 0$ for all $x\neq y$, $x,y\in X$;

\noindent(ii) there exists a constant $a_1> 0$, such that
$\rho(x,y)\leq a_1(\rho(x,z)+\rho(z,y))$ for all $x,\,y,\,z\in X$;

\noindent(iii) there exists a constant $a_0> 0$, such that
$\rho(x,y) \leq a_0 \rho(y,x) $  for all $x,\,y,\in X$.

We assume that the balls $B(x,r):=\{y\in X:\rho(x,y)<r\}$ are
measurable and $0\leq\mu (B(x,r))<\infty$ for all $x \in X$ and
$r>0$; for every neighborhood $V$ of $x\in X,$ there exists $r>0,$
such that $B(x,r)\subset V.$ Throughout the paper we also suppose
that $\mu \{x\}=0$ and that
$$  B(x,R) \setminus B(x, r) \neq \emptyset \eqno{(1.1)}$$
for all  $x\in X$, positive $r$ and $R$ with  $0< r < R < d$,
where
$$d:= diam \; (X) = \;\sup\{ \rho(x,y): x,y\in X\}. $$

Throughout the paper we suppose that $d<\infty$ and that the
doubling condition

$$  \mu(B(x, 2r))\leq c \mu (B(x,r))   $$
for $\mu$ is satisfied, where the positive constant $c$ does not
depend on $x\in X$ and $r>0$. In this case the triple  $(X, d,
\mu)$ is called a space of homogeneous type $(\hbox{SHT})$. For
the definition, examples and some properties of an $SHT$ see,
e.g., monographs \cite{StTo}, \cite{CoWe}.

A quasi-metric measure space, where the doubling condition is not
assumed is  called a non-homogeneous space.

Notice that the condition $d<\infty$ implies that $\mu(X)<\infty$
because every ball in $X$ has a finite measure.

We say that the measure $\mu$ is upper Ahlfors $Q$-- regular if
there is a positive constant $c_1$ such that $ \mu B(x,r) \leq c_1
r^Q $ for for all $x\in X$ and $r>0$. Further, $\mu$ is lower
Ahlfors $q-$ regular if there is a positive constant $c_2$ such
that $ \mu B(x,r) \geq c_2 r^q  $ for all $x\in X$ and $r>0$.

Let $1<p<\infty$, $\theta>0$ and $0\leq \lambda<1$. We denote by
$L^{p), \theta, \lambda}(X, \mu)$ the class of those  $f: X\to
{\Bbb{R}}$ for which the norm
$$
\|f\|_{L^{p), \theta, \lambda}(X, \mu)}=\sup\limits_{0<\varepsilon<p-1} \sup_{\substack{x\in X\\
0<r<d
}}\left[\frac{\varepsilon^{\theta}}{(\mu(B(x,r))^{\lambda}}\int\limits_{B(x,r)}|f(y)|^{p-\varepsilon}d\mu(y)\right]^
{\frac{1}{p-\varepsilon}}
$$
is finite.

If $\lambda =0$, then $L^{p), \theta, \lambda}(X, \mu)$ is the grand lebesgue space defined on $X$ and denoted by $L^{p), \theta}(X, \mu)$. Further,
if $\theta=1$, then we use the symbol $L^{p),  \lambda}(X, \mu)$
instead of $L^{p), \theta, \lambda}(X, \mu)$.

Using H\"older's inequality it is easy to see that the following
embeddings hold for $L^{p), \theta}$ spaces (see also
\cite{FiGuJa}, \cite{GrIwSb}):
$$ L_{w}^{p}(X, \mu)\subset L_{w}^{p),\theta_{1}}(X, \mu)\subset  L_{w}^{p),\theta_{2}}(X, \mu)\subset
L_{w}^{p-\varepsilon}(X, \mu),$$ where $0< \varepsilon< p-1$ and
$\theta_1 < \theta_2. $

The classical Morrey space, denoted by $L^{p, \lambda}(X, \mu)$,
is defined by the norm
$$
\|f\|_{L^{p, \lambda}(X, \mu)}= \sup_{\substack{x\in X\\
0<r<d
}}\left[\frac{1}{(\mu(B(x,r))^{\lambda}}\int\limits_{B(x,r)}|f(y)|^{p}d\mu(y)\right]^
{\frac{1}{p}}.
$$

Finally we mention that constants (often different constants in
the same series of inequalities) will generally be denoted by $c$
or $C$. By the symbol $p'$ we denote the conjugate number of $p$,
i.e. $p' := \frac{p}{p-1}$, $1<p<\infty$.

\section{Maximal Operator in Grand Morrey Spaces}

In this section we prove the boundedness of the Hardy--Littlewood
maximal operator

$$
(Mf)(x)=\sup_{\substack{x\in X\\ 0\leq r<d }}\frac{1}{\mu(B(x,r))}
\int\limits_{B(x,r)}|f(y)|d\mu(y), \;\; x\in X,
$$
in $L^{p),\theta, \lambda}(X,\mu)$.

Our main theorem in this section is the following statement:

{\bf Theorem 2.1.} {\em  Let $1<p<\infty$, $\theta>0$  and let
$0\leq\lambda<1$. Suppose that $d<\infty$. Then the
Hardy-Littlewood maximal operator $M$ is bounded in $L^{p),\theta,
\lambda}(X,\mu)$.}

To prove Theorem 2.1 we need some auxiliary statements.

\vskip+0.2cm

{\bf Proposition 2.1. } {\em  Let $1<p<\infty$. Then there is a
positive constant $c_{0}$ non-depending on $p$ such that
$$
\|Mf\|_{L^{p}(X,\mu)}\leq
c_{0}\left(p'\right)^{\frac{1}{p}}\|f\|_{L^{p}(X,\mu)}
\eqno{(2.1)}
$$}

\begin{proof}
The proof follows directly from the Marcinkiewicz interpolation
theorem. The constant $c_{0}$ arises from the appropriate covering
lemma (see, e. g., \cite{Du}, p. 29).
\end{proof}
Let us denote by $L^{p,\lambda}(X,\mu)$ the classical Morrey
space, where $1<p<\infty$ and $0\leq\lambda<1$, which is the class
of all $\mu$-measurable functions $f$ for which the norm
$$
\|f\|_{L^{p,\lambda}(X,\mu)}=\sup_{\substack{x\in X\\ 0\leq r<d }}
\left[\frac{1}{(\mu
B(x,r))^{\lambda}}\int\limits_{B(x,r)}|f(y)|^{p}d\mu(y)\right]^{\frac{1}{p}}
$$
is finite. If $\lambda=0$, then
$L^{p,\lambda}(X,\mu)=L^{p}(X,\mu)$.

%Proposition 3.1.2

\vskip+0.2cm

{\bf Proposition 2.2.} {\em Let $1<p<\infty$ and let
$0\leq\lambda<1$. Then
$$
\|Mf\|_{L^{p,\lambda}(X,\mu)}\leq
\left(b^{\lambda/p}c_{0}\left(p'\right)^{\frac{1}{p}}+1\right)\|f\|_{L^{p,\lambda}(X,\mu)}
$$
holds, where the positive constant $b$ arises in the doubling
condition for $\mu$ and $c_{0}$ is the constant from $(2.1)$.}

\begin{proof}
Let $r$ be $a$ small positive number and let us represent $f$ as
follows:
$$
f=f_{1}+f_{2},
$$
where $f_{1}=f\cdot\chi_{B(x,\overline{a}r)}$, $f_{2}=f-f_{1}$ and
$a$ is the positive constant given by
$\overline{a}=a_{1}(a_{1}(a_{0}+1)+1)$ (here $a_{0}$ and $a_{1}$
are constants arisen in the triangle inequality for the
quasi-metric $\rho$).

We have
$$ \left[\frac{1}{(\mu B(x,r))^{\lambda}}\int\limits_{B(x,r)}(Mf)^{p}(y)d\mu(y)\right]^{\frac{1}{p}}
\leq\left(\frac{1}{(\mu
B(x,r))^{\lambda}}\int\limits_{B(x,r)}(Mf_{1})^{p}(y)d\mu(y)\right)^{\frac{1}{p}} $$

$$ +\left(\frac{1}{(\mu B(x,r))^{\lambda}}\int\limits_{B(x,r)}(Mf_{2})^{p}(y)d\mu(y)\right)^{1/p}=:
J_{1}(x,r)+J_{2}(x,r). $$
By applying Proposition 2.1 we have that
$$ J_{1}(x,r)\leq \frac{1}{(\mu B(x,r))^{\lambda/p}}\left(\int\limits_{X}(Mf_{1}(y))^{p}d\mu(y)\right)^{1/p} $$

$$ \leq c_{0}(p')^{\frac{1}{p}}(\mu B(x,r))^{-\lambda/p}\left(\int\limits_{B(x,\overline{a}r)}|f(y)|^{p}d\mu(y)\right)^{1/p}
\leq c_{0}b^{\frac{\lambda}{p}}(p')^{\frac{1}{p}}\|f\|_{L^{p,\lambda}(X,\mu)}, $$
where  $c_{0}$ is the constant from (2.1) and $b$ arises from the
inequality
$$
\mu B(x,\overline{a}r)\leq b \mu B(x,r)
$$
which is a consequence of the doubling condition. Further, observe
that (see also \cite{KokMesAJM}, p. 23) if $y\in B(x,r)$, then
$B(x,r)\subset B(y,a_{1}(a_{0}+1)r)\subset B(x,\overline{a},r)$.
Hence, if $y\in B(x,r)$, then
$$ Mf_{2}(y)\leq\sup\limits_{B(x,r)\subset B}\frac{1}{\mu (B)}
\int\limits_{B}|f(y)|d\mu(y). $$

Consequently,

$$ J_{2}(x,r)\leq \mu (B(x,r))^{\frac{1-\lambda}{p}} \sup\limits_{B(x,r)\subset B}\left(\frac{1}{\mu (B)}
\int\limits_{B}|f(y)|^{p} d\mu(y) \right)^{1/p}
$$

$$ \leq\sup\limits_{B}(\mu B)^{-\lambda/p}\left(\int\limits_{B}|f(y)|^{p}d\mu(y)
\right)^{1/p}=\|f\|_{L^{p,\lambda}(X,\mu)}. $$
Taking into account the estimates for $J_{1}(x,r)$ and
$J_{2}(x,r)$ we conclude that

$$ \left(\frac{1}{\mu (B(x,r))^{\lambda}}\int\limits_{B(x,r)}(Mf(y))^{p}d\mu(y)\right)^{1/p} \leq\left(c_{0}b^{\lambda/p}(p')^{1/p}+1\right)\|f\|_{L^{p,\lambda}(X,\mu)}. $$
\end{proof}

\vskip+0.2cm {\em Proof of Theorem} 2.1. It is obvious that
$$
\|Mf\|_{L^{p),\theta,
\lambda}(X,\mu)}=\max\bigg\{\sup\limits_{0<\varepsilon\leq\sigma}
\sup_{\substack{x\in X\\ 0\leq r<d
}}\left(\frac{\varepsilon^{\theta}}{(\mu
B(x,r))^{\lambda}}\int\limits_{B(x,r)}(Mf(y))^{p-\varepsilon}d\mu(y) \right)^{\frac{1}{p-\varepsilon}};
$$
$$
\sup\limits_{\sigma<\varepsilon<p-1} \sup_{\substack{x\in X\\
0\leq r<d }}\left(\frac{\varepsilon^{\theta}}{(\mu
B(x,r))^{\lambda}}\int\limits_{B(x,r)}(Mf(y))^{p-\varepsilon}d\mu(y)\right)^{\frac{1}{p-\varepsilon}}\bigg\}
=:\max\left\{A_{1},A_{2}\right\}.
$$

We begin to estimate $A_{2}$. Using the facts that
$\sup\limits_{\sigma\leq\varepsilon<p-1}\varepsilon^{\frac{1}{p-\varepsilon}}=p-1$,
$\frac{1}{p-\varepsilon}>\frac{1}{p-\sigma}$ (when
$\sigma<\varepsilon<p-1$) and H\"{o}lder's inequality we have that
$$
A_{2}=\sup\limits_{\sigma<\varepsilon<p-1}\varepsilon^{\frac{\theta}{p-\varepsilon}}
\|Mf\|_{L^{p-\varepsilon,\lambda}(X,\mu)}
$$

$$
=\sup\limits_{\sigma<\varepsilon<p-1}\varepsilon^{\frac{\theta}{p-\varepsilon}}
\sup_{\substack{x\in X\\ 0\leq r<d }}(\mu
B(x,r))^{\frac{1-\lambda}{p-\varepsilon}}\left(\frac{1}{\mu
B(x,r)}\int\limits_{B(x,r)}(Mf(y))^{p-\varepsilon}d\mu(y) \right)^{\frac{1}{p-\varepsilon}}
$$

$$
\leq\sup\limits_{\sigma<\varepsilon<p-1}\varepsilon^{\frac{\theta}{p-\varepsilon}}
\sup_{\substack{x\in X\\ 0\leq r<d }}(\mu
B(x,r))^{\frac{1-\lambda}{p-\varepsilon}}\left(\frac{1}{\mu
B(x,r)}\int\limits_{B(x,r)}(Mf(y))^{p-\sigma}d\mu(y) \right)^{\frac{1}{p-\sigma}}
$$

$$
\leq\left(\sup\limits_{\sigma<\varepsilon<p-1}\varepsilon^{\frac{\theta}{p-\varepsilon}}\right)
\left(\sup_{\substack{x\in X\\ 0\leq r<d}} (\mu
B(x,r))^{\frac{1-\lambda}{p-\sigma}}\left(\frac{1}{\mu B(x,r)}
\int\limits_{B(x,r)}(M f(y))^{p-\sigma}d\mu(y) \right)^{\frac{1}{p-\sigma}}\right)
$$

$$
\leq (p-1)^{\theta}\sup_{\substack{x\in X\\ 0\leq r<d }}(\mu
B(x,r))^{\frac{1-\lambda}{p-\sigma}}\sigma^{-\frac{\theta}{p-\sigma}}\sigma^{\frac{\theta}{p-\sigma}}
\left(\frac{1}{\mu
B(x,r)}\int\limits_{B(x,r)}(Mf(y))^{p-\sigma}d\mu(y) \right)^{\frac{1}{p-\sigma}}
$$

$$
=(p-1)^{\theta}\sigma^{-\frac{\theta}{p-\sigma}}\sup_{\substack{x\in X\\
0\leq r<d }}\left(\frac{\sigma^{\theta}}{\mu
B(x,r)^{\lambda}}\int\limits_{B(x,r)}(Mf(y))^{p-\sigma}d\mu(y) \right)^{\frac{1}{p-\sigma}}
$$

$$
\leq(p-1)^{\theta}\sigma^{-\frac{\theta}{p-\sigma}}\sup\limits_{0<\varepsilon\leq\sigma}
\varepsilon^{\frac{\theta}{p-\varepsilon}}
\|Mf\|_{L^{p-\varepsilon,\lambda}(X,\mu)}.
$$
Hence, by using Proposition 2.2 we find that
$$
\|Mf\|_{L^{p),\lambda}(X,\mu)}\leq
p\sigma^{-\frac{\theta}{p-\sigma}}
\sup\limits_{0<\varepsilon\leq\sigma}\varepsilon^{\frac{\theta}{p-\varepsilon}}
\|Mf\|_{L^{p-\varepsilon,\lambda}(X,\mu)}
$$

$$
\leq
c_{0}p\cdot\sigma^{-\frac{\theta}{p-\sigma}}\sup\limits_{0<\varepsilon\leq\sigma}b^{\frac{\lambda}{p-\varepsilon}}
\left[\left(\frac{p-\varepsilon}{p-\varepsilon-1}\right)^{\frac{1}{p-\varepsilon}}+1\right]\varepsilon^{\frac{\theta}{p-\varepsilon}}
\|f\|_{L^{p-\varepsilon,\lambda}(X,\mu)}
$$

$$
\leq
c_{0}p\cdot\sigma^{-\frac{\theta}{p-\sigma}}\left[\sup\limits_{0<\varepsilon\leq\sigma}b^{\frac{\lambda}{p-\varepsilon}}
\left[\left(\frac{p-\varepsilon}{p-\varepsilon-1}\right)^{\frac{1}{p-\varepsilon}}+1\right]\right]
\|f\|_{L^{p),\lambda}(X,\mu)}.
$$

Since $\sigma$ is sufficiently small, we have that the expression

$$ S_{p,\sigma}:=c_{0}p \sigma^{-\frac{\theta}{p-\sigma}}\sup\limits_{0<\varepsilon\leq\sigma}b^{\frac{\lambda}{p-\varepsilon}}
\left[ \left( (p-\varepsilon)' \right)^{\frac{1}{p-\varepsilon}}+1\right] $$
is finite.

In fact,
$$ S_{p,\sigma}\leq c_{0}p \sigma^{-\frac{\theta}{p-\sigma}}b^{\frac{\lambda}{p-\sigma}}
\left[(p-\sigma)'+1\right]. $$

Finally,
$$
\|Mf\|_{L^{p),\theta,
\lambda}(X,\mu)}\leq\left(\inf\limits_{0<\sigma<p-1}S_{p,\sigma}\right)
\|f\|_{L^{p),\theta, \lambda}(X,\mu)}.
$$
$\Box$

\

\section{Calder\'{o}n-Zygmund Operators in Grand Morrey Spaces}

Let
$$ Tf (x)= p.v. \int\limits_X k(x,y)f(y) d\mu(y), $$
where  $k: X\times X\setminus \{(x,x): x\in X\} \to {\Bbb{R}}$ be
a measurable function satisfying the conditions:
\begin{gather*}
|k(x,y)|\leq \frac{c}{\mu B(x, \rho(x,y))}, \;\; x,y\in X, \;\; x\neq y; \\
 |k(x_1,y)-k(x_2,y)|+ |k(y, x_1)-k(y, x_2)| \leq c \omega \Big( \frac{\rho(x_2, x_1)}{\rho(x_2, y)}\Big)
 \frac{1}{\mu B(x_2, \rho(x_2,y))}
 \end{gather*}
for all $x_1, x_2$ and $y$ with $\rho(x_2,y)> c \rho(x, x_2)$,
where $\omega$ is a positive non-decreasing function on
$(0,\infty)$ which satisfies the  $\Delta_2$ condition:
$\omega(2t)\leq c \omega(t)$ ($t>0$);  and the Dini condition:
$\int_0^1 \big(\omega(t)/t\big) dt <\infty$.

We also assume that for some constant $p_0$, $1<p_0<\infty$, and
all $f\in L^{p_0}(X, \mu)$ the limit $Tf(x)$ exists almost
everywhere on $X$ and that $T$ is bounded in $L^{p_0}(X, \mu)$.

\vskip+0.1cm

For simplicity we will assume that $p_{0}=2$. We call $T$ the
Calder\'{o}n-Zygmund operator.

\vskip+0.2cm

{\bf Lemma 3.1.} {\em Let $T$ be the Calder\'{o}n-Zygmund
operator. Then there is a positive constant $c$ non-depending on
$p$ such that the following estimates hold:
$$
\|T\|_{L^{p}(X,\mu)\rightarrow L^{p}(X,\mu)}\leq
c\left(\frac{p}{p-1}+\frac{p}{2-p}\right),\;\;1<p<2,
$$
$$
\|T\|_{L^{p}(X,\mu)\rightarrow L^{p}(X,\mu)}\leq
c\left(p+\frac{p}{p-2}\right),\;\;p>2.
$$}

\begin{proof}
Since $T$ has weak $(1,1)$ and strong $(2,2)$ types the
Marcinkiwicz interpolation theorem (see, e. g., \cite{Du}, p. 29)
we have that
$$
\|T\|_{L^{p}(X,\mu)\rightarrow L^{p}(X,\mu)}\leq
\left(\frac{2p}{p-1}\frac{A_{0}}{c^{p-1}}+\frac{4p}{2-p}
\frac{A_{1}^{2}}{c^{p-2}}\right)^{\frac{1}{p}}\|f\|_{L^{p}(X,\mu)},\;\;1<p<2,
$$
where $A_{0}$ is the constant arisen in the weak $(1,1)$ type
inequality for $T$ and $A_{1}$ is the constant from the strong
$(2,2)$ type inequality for $T$. Observe now that
$$
\left[\frac{2p}{p-1}\frac{A_{0}}{c^{p-1}}+\frac{4p}{2p}
\frac{A_{1}^{2}}{c^{p-2}}\right]^{1/p} \leq
2^{1/p}\left(\frac{p}{p-1}\right)^{1/p}\frac{A_{0}^{1/p}}{c^{(p-1)/p}}+4^{1/p}
\left(\frac{p}{2-p}\right)^{1/p}\frac{A_{1}^{2/p}}{c^{(p-2)/p}}
$$
$$
\leq c \left(\frac{p}{p-1}+\frac{p}{2-p}\right),
$$
where the positive constant $c$ does not depend on $p$.

Let now $p>2$. By using the above--mentioned arguments we have that

$$ \|T\|_{L^{p}\rightarrow L^{p}}=\|T\|_{L^{p'}\rightarrow L^{p'}}\leq \left(\frac{p'}{p'-1}+\frac{p'}{2-p'}\right)=c\left(p+\frac{p}{p-2}\right).
$$
\end{proof}

%Proposition 3.2.1
{\bf Proposition 3.1.} {\em  Let $1<p<\infty$, $0\leq\lambda<1$.
Then
$$
\|T\|_{L^{p,\lambda}(X,\mu)}\leq
c\left[\frac{p}{p-1}+\frac{p}{2-p}+\frac{p-\lambda+1}{1-\lambda}\right],\;\;\;1<p<2,
$$
$$
\|T\|_{L^{p,\lambda}(X,\mu)}\leq
c\left[p+\frac{p}{p-2}+\frac{p-\lambda+1}{1-\lambda}\right],
\;\;\;p>2.
$$
}

\begin{proof}
Let us take small $r>0$ and $x\in X$. Represent $f$ as follows:
$f=f_{1}+f_{2}$, where $f_{1}=f\cdot\chi_{B(x,2a_{1}r)}$,
$f_{2}=f-f_{1}$, where $a_{1}$ is the constant from the triangle
inequality for the quasi-metric $\rho$. Observe that if $y\in
B(x,r)$ and $z\in X\backslash B(x,2a_{1}r)$, then
$$ \mu B(x,\rho(x,z))\leq c\mu B(y,\rho(y,z)). \eqno{(3.1)} $$
Inequality (3.1) follows from the estimates

$$ \mu B(x,\rho(x,z))\leq c_{1}\mu B(x,\rho(y,z))\leq c_{2}\mu B(y,\rho(y,z)). \eqno{(3.2)} $$
To show the first part of (3.2) observe that
$$\rho(x,z)\leq a_{1}\rho(x,y)+a_{1}\rho(y,z)\leq a_{1}r+a_{1}\rho(y,z)
\leq\frac{\rho(x,z)}{2}+a_{1}\rho(y,z)
$$
Hence, $\frac{\rho(x,z)}{2a_{1}}\leq\rho(y,z)$. Now by the doubling condition we have the first part of (3.2).

The second part of (3.2) follows easily.

Recall now that the doubling condition for $\mu$ implies the
reverse doubling condition for $\mu$: there one constants
$0<\alpha,\beta<1$ such that for all $x\in X$ and small positive
$r$,
$$
\mu B(x,\alpha r)\leq\beta\mu B(x,r).\;\;\;\;\;\;\eqno{(3.3)}
$$
Let us take an integer $m_{0}$ so that $\alpha^{m_{0}}d$ is
sufficiently small, where $d$ is the diameter of $X$.

Let $y\in B(x,r)$. Then by (3.1) and Fubini's theorem we have that
$$ |Tf_{2}(y)|\leq c\int\limits_{X\backslash B(x,2a_{1}r)} |f(z)|(\mu
B(x,z))^{-1}d\mu(z) $$

$$ \leq c_{\beta}\int\limits_{X\backslash B(x,2a_{1}r)}|f(z)|\left(\int\limits_{B(x,\alpha^{m_{0}}\rho(x,z))\backslash
B(x,\alpha^{m_{0}-1}\rho(x,z))}(\mu B(x,\rho(x,t))^{-2}d\mu(t)\right)d\mu(z) $$

$$ \leq c_{\beta}\int\limits_{X\backslash B(x,2\alpha^{m_{0}-1}a_{1}r)}(\mu
B(x,\rho(x,t)))^{-2} \left(\int\limits_{B(x,\alpha^{1-m_{0}}\rho(x,t))} |f(z)|d\mu(z)\right)d\mu(t) $$

$$ \leq c \int\limits_{X\backslash B(x,2\alpha^{m_{0}-1}a_{1}r)}(\mu
B(x,\rho(x,t)))^{-1}\overline{f}(x,t)d\mu(t), $$
where

$$ \overline{f}(x,t):=\left[\mu B\left(x,\alpha^{1-m_{0}}r\right)\right]^{-1}
\int\limits_{B(x,\alpha^{1-m_{0}}\rho(x,t))}|f(z)|d\mu(z). $$
Observe that by H\"{o}lder's inequality the following estimates
hold:

$$ \overline{f}(x,t)\leq\mu B(x,\alpha^{1-m_{0}}\rho(x,t))^{-1} \|f\|_{L^{p}(B(x,\alpha^{1-m_{0}}\rho(x,t))} \|\chi
_{B(x,\alpha^{1-m_{0}}\rho(x,t))}\|_{L^{p'}(X)} $$

$$ \leq c (\mu B(x,\alpha^{1-m_{0}}\rho(x,t))^{-\frac{1}{p}}
\|f\|_{L^{p}(B(x,\alpha^{1-m_{0}}\rho(x,t))} $$

$$ \leq c (\mu B(x,\alpha^{1-m_{0}}\rho(x,z))^{\frac{\lambda-1}{p}}
\|f\|_{L^{p,\lambda}(X,\mu)}. $$ By applying now Lemma 1.2 of
\cite{KokMesAJM} (see also the monograph \cite{EdKoMe}. p. 372) we
find that for $y\in B(x,r)$,

$$ |Tf_{2}(y)|\leq c\|f\|_{L^{p,\lambda}(X,\mu)}
\int\limits_{X\backslash B(x,2\alpha^{1-m_{0}}a_{1}r)} \left[\mu
B(x,\alpha^{1-m_{0}}\rho(x,t))\right]^{\frac{\lambda-1}{p}-1}d\mu(t)
$$

$$
\leq c\left[\mu B(x,2
\alpha^{1-m_{0}}a_{1}r)\right]^{\frac{\lambda-1}{p}}\cdot\frac{\frac{\lambda-1}{p}-1}{\frac{\lambda-1}{p}}
\cdot\|f\|_{L^{p,\lambda}(X,\mu)} $$

$$ =c\frac{p-\lambda+1}{1-\lambda}\left[\mu
B(x,2\alpha^{1-m_{0}}a_{1}r)\right]^{\frac{\lambda-1}{p}} \| f
\|_{L^{p,\lambda}(X, \mu)},
$$ where the positive constant $c$ does not depend on $\lambda$
and $p$.

Consequently, by Lemma 3.1 we find that
$$ \left[(\mu B(x,r))^{-\lambda}\int\limits_{B(x,r)}|Tf(y)|^{p}d\mu(y)\right]^{1/p}
\leq (\mu B(x,r))^{-\lambda/p}\left(\int\limits_{B(x,r)}|Tf_{1}(y)|^{p}d\mu(y)\right)^{\frac{1}{p}}
$$

$$ +(\mu B(x,r))^{-\lambda/p}\left(\int\limits_{B(x,r)}|Tf_{2}(y)|^{p}d\mu(y)\right)^{\frac{1}{p}}$$

$$ \leq c\cdot c_{p}(\mu B(x,r))^{-\frac{\lambda}{p}} \left(\int\limits_{B(x,2a_{1}r)}|f(y)|^{p}d\mu(y)\right)^{\frac{1}{p}} $$

$$ +c\frac{p-\lambda+1}{1-\lambda}\left(\mu B(x,2\alpha^{1-m_{0}}a_{1}r)\right)^{\frac{\lambda-1}{p}} (\mu B(x,2a_{1}r))^{\frac{1}{p}} $$

$$ \times(\mu B(x,r))^{-\frac{\lambda}{p}}\|f\|_{L^{p,\lambda}(X,\mu)} $$

$$ \leq\left(c\cdot c_{p}+c\cdot\frac{p-\lambda+1}{1-\lambda}\right)\|f\|_{L^{p,\lambda}(X,\mu)} =
c\left(c_{p}+\frac{p-\lambda+1}{1-\lambda}\right)\|f\|_{L^{p,\lambda}(X,\mu)}, $$
where the constant $c>0$ does not depend on $p$ and $\lambda$, and
$$
c_{p}=\left\{%
\begin{array}{ll}
    \frac{p}{p-1}+\frac{p}{2-p},\;\;\;\;1<p<2, \\
    p+\frac{p}{p-2},\;\;\;\;\;\;\;p>2. \\
\end{array}%
\right.
$$
Observe that the constant $c$ does not depend on $\lambda$ and
$p$.
\end{proof}
\vskip+0.1cm

{\bf Theorem 3.1.} {\em Let $1<p<\infty$, $\theta>0$ and let
$0<\lambda<1$. Then the operator $T$ is bounded in $L^{p),\theta,
\lambda}(X,\mu)$. }

\begin{proof}
We have
$$
\|Tf\|_{L^{p),\theta, \lambda}(X,\mu)}=\max\bigg\{\sup\limits_{0<\varepsilon\leq\sigma}
\sup_{\substack{x\in X\\ 0\leq
r<d}}\left(\frac{\varepsilon^{\theta}}{(\mu
B(x,r))^{\lambda}}\int\limits_{B(x,r)}|Tf|^{p-\varepsilon}\right)^{\frac{1}{p-\varepsilon}},
$$

$$
\sup\limits_{\sigma<\varepsilon<p-1} \sup_{\substack{x\in X\\
0\leq r<d}}\left(\frac{\varepsilon^{\theta}}{(\mu
B(x,r))^{\lambda}}\int\limits_{B(x,r)}|Tf|^{p-\varepsilon}\right)^{\frac{1}{p-\varepsilon}}\bigg\}
=:\max\{A_{1},A_{2}\},
$$
where $\sigma$ is the number satisfying the condition
$0<\sigma<p-1$.

Now we estimate $A_{2}$. By applying H\"{o}lder's inequality we find that
$$
A_{2}=\sup\limits_{\sigma<\varepsilon<p-1}\varepsilon^{\frac{\theta}{p-\varepsilon}}
\|Tf\|_{L^{p-\varepsilon,\lambda}(X,\mu)}
=\sup\limits_{\sigma<\varepsilon<p-1}\varepsilon^{\frac{\theta}{p-\varepsilon}}
\mu (B(x,r))^{-\frac{\lambda}{p-\varepsilon}}(\mu
B(x,r))^{\frac{1}{p-\varepsilon}} $$

$$ \times\left(\frac{1}{\mu B(x,r)}\int\limits_{B(x,r)}|Tf|^{p-\varepsilon}\right)^{\frac{1}{p-\varepsilon}} $$

$$ \leq\sup\limits_{\sigma<\varepsilon<p-1}\varepsilon^{\frac{\theta}{p-\varepsilon}}(\mu B(x,r))^{\frac{1-\lambda}{p-\varepsilon}}\left( (\mu
B(x,r))^{-1}\int\limits_{B(x,r)}|Tf|^{p-\sigma}\right)^{\frac{1}{p-\sigma}}. $$

Further, without loss of generality we can assume that $\mu(X)=1$ and, consequently,
$(\mu(B(x,r))^{\frac{1-\lambda}{p-\varepsilon}}\leq (\mu
B(x,r))^{\frac{1-\lambda}{p-\sigma}}$. Hence,

$$ A_2 \leq (p-1)^{\theta}\sigma^{-\frac{\theta}{p-\sigma}}\sup\limits_{0<\varepsilon\leq\sigma}
\varepsilon^{\frac{\theta}{p-\varepsilon}}
\|Tf\|_{L^{p-\varepsilon,\lambda}(X,\mu)}. $$

Hence, by Proposition 3.1 we conclude that

$$ \|Tf\|_{L^{p),\theta, \lambda}(X,\mu)}\leq\left[(p-1)^{\theta}\sigma^{-\frac{\theta}{p-\sigma}}+1\right]
\sup\limits_{0<\varepsilon\leq\sigma}
\varepsilon^{\frac{\theta}{p-\varepsilon}}
\|Tf\|_{L^{p-\varepsilon,\lambda}(X,\mu)} $$

$$ \leq\left[(p-1)^{\theta}\sigma^{-\frac{\theta}{p-\sigma}}+1\right]\sup\limits_{0<\varepsilon\leq\sigma}
C_{p,\lambda,\varepsilon}\varepsilon^{\frac{\theta}{p-\varepsilon}}
\|f\|_{L^{p-\varepsilon,\lambda}(X,\mu)} $$

$$ =\left[(p-1)^{\theta}\sigma^{-\frac{\theta}{p-\sigma}}+1\right]\sup\limits_{0<\varepsilon\leq\sigma}
C_{p,\lambda,\varepsilon} \sup_{\substack{x\in X\\ 0\leq
r<d}}\left[\frac{\varepsilon^{\theta}}{(\mu
B(x,r))^{\lambda}}\int\limits_{B(x,r)}|f(y)|^{p-\varepsilon}d\mu(y)  \right]^{\frac{1}{p-\varepsilon}}
$$

$$
\leq\left[(p-1)^{\theta}\sigma^{-\frac{\theta}{p-\sigma}}+1\right]\|f\|_{L^{p),\theta,
\lambda}(X,\mu)} \sup\limits_{0<\varepsilon\leq\sigma}
C_{p,\lambda,\varepsilon}, $$
where

$$ C_{p,\lambda,\varepsilon}=\left\{%
\begin{array}{ll}
    \frac{p-\varepsilon-\lambda+1}{1-\lambda}+\frac{p-\varepsilon}{p-\varepsilon-1}+\frac{p-\varepsilon}{2-p+\varepsilon}, \;\;\;\;1<p<2 \\
    \frac{p-\varepsilon-\lambda+1}{1-\lambda}+p-\varepsilon+\frac{p-\varepsilon}{p-\varepsilon-2},\;\;\;\;\;p>2. \\
\end{array}%
\right. $$

Observe now that
$$
\sup\limits_{0<\varepsilon\leq\sigma}C_{p,\lambda,\varepsilon}\leq
\left\{%
\begin{array}{ll}
    \frac{p-\lambda+1}{1-\lambda}+\frac{p-\sigma}{p-\sigma-1}+\frac{p}{2-p},\;\;\;\;\;\;\;1<p<2, \\
    \frac{p-\lambda+1}{1-\lambda}+\frac{p-\sigma}{p-\sigma-1}+\frac{p-\sigma}{p-\sigma-2},\;\;\;\;p>2, \\
\end{array}%
\right.
$$
where $\sigma$ is sufficiently small. \
\end{proof}

\section{Fractional integrals  in grand Morrey spaces}

\subsection{ Potentials $(I_{\alpha}f)(x)=\int\limits_{X}\frac{f(y)}{\rho(x,y)^{\gamma-\alpha}}d\mu(y)$}

Let an SHT $(X,\rho,\mu)$ satisfy the condition: there are
positive constants $b$ and  $\gamma$ such that
$$
\mu B(x,r)\leq b r^{\gamma}, \eqno{(4.1)}
$$
for all $x\in X$ and $r$, $0<r<d$, i.e. $\mu$ is upper $\gamma-$
Ahlfors regular. As before we assume that $d= diam(X)< \infty$.

Let
$$ (I_{\alpha}f)(x)=\int\limits_{X}\frac{f(y)}{\rho(x,y)^{\gamma-\alpha}}d\mu(y), \;\; x\in X, $$
where $0<\alpha<\gamma$.

In this section we study the boundedness of $I_{\alpha}$ in grand
Morrey spaces. For this we define the classical Morrey space as
follows: $f\in {\mathcal{L}}^{p,\lambda}(X,\mu)$ ($1<p<\infty$,
$0\leq\lambda<\frac{1}{\gamma}$) if
$$
\|f\|_{{\mathcal{L}}^{p,\lambda}(X,\mu)}:=\sup_{\substack{x\in X\\ 0\leq r<d}}
\left(\frac{1}{r^{\gamma\lambda}}\int\limits_{B(x,r)}|f(y)|^{p}d\mu(y)\right)^{1/p}<\infty
$$
Further,  let $\varphi$ be a positive function in $(0, p-1)$ which
is increasing near $0$ and satisfies the condition $\varphi
(0+)=0$. We say that $f\in
{\mathcal{L}}^{p),\varphi(\cdot),\lambda}(X,\mu)$ if
$$
\|f\|_{{\mathcal{L}}^{p),\varphi(\cdot),\lambda}(X,\mu)}=\sup\limits_{0<\varepsilon<p-1}\sup_{\substack{x\in
X\\ 0\leq
r<d}}\left(\frac{\varphi(\varepsilon)}{r^{\gamma\lambda}}\int\limits_{B(x,r)}|f|^{p}d\mu\right)^{\frac{1}{p}}<\infty.
$$

Let $\theta$ be a positive number. If  $\varphi(\varepsilon)\equiv
\varepsilon^{\theta}$, then we denote
${\mathcal{L}}^{p),\varphi(\cdot),\lambda}(X,\mu)=:{\mathcal{L}}^{p),\theta,
\lambda}(X,\mu)$. For $\theta=1$ we have the grand Morrey space
${\mathcal{L}}^{p), \lambda}(X, \mu)$. \vskip+0.2cm

Let
$$ ({\mathcal{M}}f)(x)=\sup_{\substack{x\in X\\ 0\leq r<d}}\frac{1}{r^{\gamma}}\int\limits_{B(x,r)}|f(y)|d\mu(y), \;\; x\in X. $$
We begin with the following statement:
\vskip+0.2cm

{\bf Proposition 4.1.} {\em Let $1<p<\infty$ and let $0\leq\lambda<1$. Then

$$ \|{\mathcal{M}}f\|_{L^{p,\lambda}(X,\mu)}\leq \left( (\overline{a})^{\frac{\lambda\gamma}{p}}c_{0}
\left(p'\right)^{\frac{1}{p}}+1\right)
\|f\|_{L^{p,\lambda}(X,\mu)} $$ holds, where  $c_{0}$ is the
constant from $(2.1)$ and $\overline{a}=a_{1}(a_{1}(a_{0}+1)+1)$.}

\vskip+0.1cm

{\em Proof.} Since ${\mathcal{M}}f(x) \leq Mf(x)$, by Proposition
2.1 we have that
$$ \|{\mathcal{M}}\|_{L^p(X, \mu) \to L^p(X, \mu)} \leq c_0
(p')^{1/p}. $$ Repeating the proof of Proposition 2.2 we have the
desired result.  $\Box$

\vskip+0.2cm

%Lemma 3.3.1.
{\bf Lemma 4.1.} {\em  Let $1<p<\infty$,
$0<\alpha<\frac{(1-\lambda)\gamma}{p}$,
$\frac{1}{p}-\frac{1}{q}=\frac{\alpha}{(1-\lambda)\gamma}$, where
$0\leq\lambda<1/\gamma$. Then the inequality
$$
\|I_{\alpha}f\|_{{\mathcal{L}}^{q,\lambda}(X,\mu)}\leq
c(p,\alpha,\lambda,\gamma)\|f\|_{{\mathcal{L}}^{p,\lambda}(X,\mu)}
$$
holds, where the positive constant $c(p,\alpha,\lambda,\gamma)$ is
given by
$$
c(p,\alpha,\lambda,\gamma)=c\frac{(1-\lambda)\gamma}{\alpha[(1-\lambda)\gamma-\alpha
p]}\left[(p')^{1/q}+1\right],
$$
and the positive constant $c$ does not depend on $p$  and
$\alpha$. }

\begin{proof}
First we show that the Hedberg's \cite{Hed} type inequality holds:
$$
|(I_{\alpha}f)(x)|\leq
c_{p,\lambda,\gamma,\alpha}({\mathcal{M}}f)^{1-\frac{p\alpha}{(1-\lambda)\gamma}}(x)\|f\|_{{\mathcal{L}}^{p,\lambda}(X,\mu)}^{\frac{\alpha
p}{(1-\lambda)\gamma}}, \eqno{(4.2)}
$$
where
$c_{p,\lambda,\gamma,\alpha}=\frac{2(1-\lambda)\gamma}{\alpha((1-\lambda)\gamma-\alpha
p))}$.
To prove (4.2) we set
$$ f_{r}(x):=\frac{1}{r^{\gamma}}\int\limits_{B(x,r)}|f(y)|d\mu(y). $$
The inequality
$$
\rho(x,y)\leq
2\;\;\;\;\;\int\limits_{\rho(x,y)}^{2\rho(x,y)}t^{\alpha-\gamma-1}dt,\;\;\;\;0<\rho(x,y)<l,
\eqno{(4.3)}
$$
is obvious. By using (4.3) we find that
$$
|(I_{\alpha}f)(x)|\leq
2\int\limits_{X}|f(y)|\left(\int\limits_{\rho(x,y)}^{2\rho(x,y)}t^{\alpha-\gamma-1}dt\right)d\mu(y)
$$
$$
=2\int\limits_{0}^{2d}t^{\alpha-\gamma-1}\left(\int\limits_{\frac{t}{2}<\rho(x,y)<t}|f(y)|d\mu(y)\right)dt
\leq 2\int\limits_{0}^{2d}t^{\alpha-\gamma}f_{t}(x)dt.
$$
Taking $\varepsilon>0$ (which will be chosen later) we have that
$$
|(I_{\alpha}f)(x)|\leq
2\left[\int\limits_{0}^{\varepsilon}t^{\alpha-\gamma}f_{t}(x)dt+\int\limits_{\varepsilon}^{2d}t^{\alpha-\gamma}f_{t}(x)dt\right]
$$
$$
=:2\left[J_{1}^{(\varepsilon)}(x)+J_{2}^{(\varepsilon)}(x)\right].
$$
It is obvious that
$$
J_{1}^{(\varepsilon)}(x)\leq
({\mathcal{M}}f)(x)\int\limits_{0}^{\varepsilon}t^{\alpha-1}dt=\frac{({\mathcal{M}}f)(x)}{\alpha}\varepsilon^{\alpha}.
$$
Further, by H\"{o}lder's inequality and condition (4.1) it is
clear that
$$
f_{t}(x)=\frac{1}{t^{\gamma}}\int\limits_{B(x,t)}|f(y)|d\mu(y)\leq\left(\frac{1}{t^{\gamma}}\int\limits_{B(x,t)}|f(y)|^{p}d\mu(y)\right)^{\frac{1}{p}}
$$

$$
=t^{-\frac{\gamma}{p}+\frac{\lambda\gamma}{p}}\left(\frac{1}{t^{\gamma\lambda}}\int\limits_{B(x,r)}|f(y)|^{p}d\mu(y)\right)^{\frac{1}{p}}
\leq
t^{-\frac{\gamma}{p}+\frac{\lambda\gamma}{p}}\|f\|_{{\mathcal{L}}^{p,\lambda}(X,\mu)}.
$$
By applying the condition $\frac{\lambda-1}{p}\gamma+\alpha<0$ we find that

$$ |(I_{\alpha}f)(x)|\leq
2\left[\frac{({\mathcal{M}}f)(x)}{\alpha}\varepsilon^{\alpha}+\left(\int\limits_{\varepsilon}^{2l}
t^{\frac{(\lambda-1)\gamma}{p}+\alpha-1}dt\right)\|f\|_{{\mathcal{L}}^{p,\lambda}}\right] $$

$$
=2\left[\frac{({\mathcal{M}}f)(x)}{\alpha}\varepsilon^{\alpha}-\frac{\varepsilon^{\frac{\lambda-1}{p}\gamma+\alpha}}{\left[\alpha+\frac{\lambda-1}{p}\gamma\right]}
\|f\|_{{\mathcal{L}}^{p,\lambda}(X,\mu)}\right].
$$
Let
$\varepsilon=\left[\frac{\|f\|_{{\mathcal{L}}^{p,\lambda}(X,\mu)}}{({\mathcal{M}}f)(x)}\right]^{\frac{p}{(1-\lambda)\gamma}}$. Then

$$ |(I_{\alpha}f)(x)|\leq 2\bigg[\frac{({\mathcal{M}}f)^{1-\frac{p\alpha}{(1-\lambda)\gamma}}(x)}{\alpha}\|f\|_
{{\mathcal{L}}^{p,\lambda}(X,\mu)}^{\frac{p\alpha}{(1-\lambda)\gamma}}
$$

$$
-\frac{1}{\left[\alpha+\frac{(\lambda-1)\gamma}{p}\right]}
\|f\|_{L^{p,\lambda}(X,\mu)}^{\frac{\alpha
p}{(1-\lambda)\gamma}}({\mathcal{M}}f)^{1-\frac{p\alpha}{(1-\lambda)\gamma}}(x)\bigg] $$

$$ =2\left[\frac{1}{\alpha}-\frac{p}{\alpha p+(\lambda-1)\gamma}\right]\|f\|_{{\mathcal{L}}^{p,\lambda}(X,\mu)}^{\frac{\alpha
p}{(1-\lambda)\gamma}}({\mathcal{M}}f)^{1-\frac{p\alpha}{(1-\lambda)\gamma}}(x) $$

$$ =2\frac{(1-\lambda)\gamma}{\alpha((1-\lambda)\gamma-\alpha p)} \|f\|_{{\mathcal{L}}^{p,\lambda}(X,\mu)}^{\frac{\alpha
p}{(1-\lambda)\gamma}}({\mathcal{M}}f)^{1-\frac{p\alpha}{(1-\lambda)\gamma}}(x). $$

Consequently, by the condition
$\frac{1}{p}-\frac{1}{q}=\frac{\alpha}{(1-\lambda)\gamma}$ and Proposition 4.1 we have
that

$$
\left(\frac{1}{t^{\gamma\lambda}}\int\limits_{B(x,t)}|(I_{\alpha}f)(y)|^{q}d\mu(y)\right)^{\frac{1}{q}}
$$

$$
\leq
t^{-\frac{\gamma\lambda}{q}}\frac{2(1-\lambda)\gamma}{\alpha((1-\lambda)\gamma-\alpha
p)}\left[\int\limits_{B(x,t)}({\mathcal{M}}f(y))^{q\left[1-\frac{p\alpha}{(1-\lambda)\gamma}\right]}d\mu(y)\right]^{\frac{1}{q}}
\|f\|_{{\mathcal{L}}^{p,\lambda}(X,\mu)}^{\frac{\alpha p}{(1-\lambda)\gamma}}
$$

$$
=\frac{2(1-\lambda)\gamma}{\alpha((1-\lambda)\gamma-\alpha
p)}\left[\frac{1}{t^{\gamma\lambda}}\int\limits_{B(x,t)}({\mathcal{M}}f(y))^{p}d\mu(y)\right]^{\frac{1}{q}}
\|f\|_{{\mathcal{L}}^{p,\lambda}(X,\mu)}^{\frac{\alpha p}{(1-\lambda)\gamma}}
$$

$$
\leq\frac{2(1-\lambda)\gamma}{\alpha\left[(1-\lambda)\gamma-\alpha
p\right]}\|{\mathcal{M}}f\|_{{\mathcal{L}}^{p,\lambda}(X,\mu)}^{p/q}\|f\|_{{\mathcal{L}}^{p,\lambda}(X,\mu)}^{\frac{\alpha
p}{(1-\lambda)\gamma}}.
$$

$$
\big(\; \text{ recall that}\;\;
\|{\mathcal{M}}f\|_{{\mathcal{L}}^{p,\lambda}(X,\mu)}\leq \left(
(\overline{a})^{\frac{\lambda\gamma}{p}}c_{0}
\left(p'\right)^{\frac{1}{p}}+1\right)
\|f\|_{{\mathcal{L}}^{p,\lambda}(X,\mu)}\big)
$$

$$
\leq \frac{2(1-\lambda)\gamma}{\alpha[(1-\lambda)\gamma-\alpha p]}
\left( (\overline{a})^{\frac{\lambda\gamma}{p}}c_{0}
\left(p'\right)^{\frac{1}{p}}+1\right)^{p/q}\|f\|_{{\mathcal{L}}^{p,\lambda}(X,\mu)}
$$

$$ \leq \frac{(1-\lambda)\gamma}{\alpha[ (1-\lambda)\gamma- \alpha p]} \left( (c_0)^{p/q}
(\overline{a})^{\frac{\lambda \gamma}{q}} (p')^{1/q} +1 \right)
\|f\|_{ {\mathcal{L}}^{p, \lambda}(X,\mu) } $$

$$ \leq \frac{(1-\lambda)\gamma}{\alpha[ (1-\lambda)\gamma- \alpha p]} \left( c_0
(\overline{a})^{\lambda \gamma} (p')^{1/q} +1 \right) \|f\|_{
{\mathcal{L}}^{p, \lambda}(X,\mu) } $$

$$ \leq c \frac{(1-\lambda)\gamma}{\alpha[ (1-\lambda)\gamma- \alpha p]} \left(  (p')^{1/q} +1 \right) \|f\|_{
{\mathcal{L}}^{p, \lambda}(X,\mu) }. $$
\end{proof}

%Theorem 4.1
{\bf Theorem 4.1.} {\em  Let Let $1<p<\infty$,
$0<\alpha<\frac{(1-\lambda)\gamma}{p}$, $0\leq\lambda<1/\gamma$
and let
$\frac{1}{p}-\frac{1}{q}=\frac{\alpha}{(1-\lambda)\gamma}$.
Suppose that $\theta_1>0$. We set
$$\theta_2= \Big[ 1+ \frac{\alpha q}{(1-\lambda)\gamma}\bigg] \theta_1.$$

Then the operator $I_{\alpha}$ is bounded from ${\mathcal{L}}^{p),\theta_1,
\lambda}(X,\mu)$ to ${\mathcal{L}}^{q),\theta_2,\lambda}(X,\mu)$. }

\begin{proof}
Let us introduce the function:
$$\varphi(u):=\left[p+\frac{(1-\lambda)(u-q)\gamma}{(1-\lambda)\gamma-\alpha(u-q)}\right]^
{\frac{\gamma(1-\lambda)-(u-q)\alpha}{(1-\lambda)\gamma}}.$$
Observe that
$$ \varphi (t) \sim t^{1+ \frac{\alpha
q}{(1-\lambda)\gamma}}, \;\;\; \text{as}\;\; t\to 0+. $$
Hence it is enough to prove that $I_{\alpha}$ is bounded from
${\mathcal{L}}^{p), \theta_1, \lambda}(X, \mu)$ to  ${\mathcal{L}}^{q), \psi(\cdot), \lambda}(X,
\mu)$, where
$$ \psi (t):= \varphi \big(t^{\theta_1}\big). $$

Let $\sigma$ be a small positive number. As in the proofs of the
main theorems of previous sections we have

$$
\|I_{\alpha}f\|_{{\mathcal{L}}^{q,\varphi,\lambda}(X,\mu)}=\max
\bigg\{\sup\limits_{0<\varepsilon\leq\sigma} \sup_{\substack{x\in
X\\ 0\leq
r<d}}\bigg(\frac{\psi(\varepsilon)}{t^{\lambda\gamma}}\int\limits_{B(x,t)}|I_{\alpha}f(x)|^{q-\varepsilon}d\mu(x)
\bigg)^{\frac{1}{q-\varepsilon}},
$$

$$ \sup\limits_{\sigma<\varepsilon<q-1} \sup_{\substack{x\in X\\
0\leq
r<d}}\bigg(\frac{\psi(\varepsilon)}{t^{\lambda\gamma}}\int\limits_{B(x,r)}|I_{\alpha}f(y)|^{q-\varepsilon}d\mu(y)\bigg)^
{\frac{1}{q-\varepsilon}}\bigg\}
=:\max\{A_{1},A_{2}\}. $$

For $A_{2}$, we observe that

$$
\psi(\varepsilon)^{\frac{1}{q-\varepsilon}}t^{\frac{(1-\lambda)\gamma}{q-\varepsilon}}
\left(\frac{1}{t^{\gamma}}\int\limits_{B(x,r)}|I_{\alpha}f(y)|^{q-\varepsilon}d\mu(y)\right)^{\frac{1}{q-\varepsilon}}
$$

$$ \leq\sup\limits_{\sigma<\varepsilon<q-1}(\psi(\varepsilon))^{\frac{1}{q-\varepsilon}}
t^{\frac{(1-\lambda)\gamma}{q-\varepsilon}}
\left(\frac{1}{t^{\gamma}}\int\limits_{B(x,r)}|I_{\alpha}f(y)|^{q-\varepsilon}d\mu(y)\right)
^{\frac{1}{q-\varepsilon}} $$

$$ \leq (\text{by H\"{o}lder's inequality and the fact that}\;\sigma<\varepsilon)$$

$$ \leq\left[\sup\limits_{\sigma<\varepsilon<q-1}(\psi(\varepsilon))^{\frac{1}{q-\varepsilon}}\right]
t^{\frac{(1-\lambda)\gamma}{q-\sigma}}t^{-\frac{\gamma}{q-\sigma}}
\left(\int\limits_{B(x,t)}|I_{\alpha}f(y)|^{q-\sigma}d\mu(y)\right)^{\frac{1}{q-\sigma}}$$

$$
=\left[\sup\limits_{\sigma<\varepsilon<q-1}(\psi(\varepsilon))^{\frac{1}{q-\varepsilon}}\right]
\psi(\sigma)^{-\frac{1}{q-\sigma}}\left(\frac{\varphi(\sigma)}{t^{\lambda\gamma}}
\int\limits_{B(x,t)}|I_{\alpha}f(y)|^{q-\sigma}d\mu(y)\right)^{\frac{1}{q-\sigma}}
$$

$$
\leq\left[\sup\limits_{\sigma<\varepsilon<q-1}\psi(\varepsilon)^{\frac{1}{q-\varepsilon}}\right]
\psi(\sigma)^{-\frac{1}{q-\sigma}}\sup\limits_{0<\varepsilon\leq\sigma} \sup_{\substack{x\in X\\
0\leq t<d}}\left(\frac{\psi(\varepsilon)}{t^{\lambda\gamma}}
\int\limits_{B(x,t)}|I_{\alpha}f(y)|^{q-\sigma}d\mu(y)\right)^{\frac{1}{q-\sigma}}. $$

Further, applying Lemma 4.1, for $\varepsilon$ satisfying the condition
$0<\varepsilon\leq\sigma$, we have

$$
\left(\frac{\psi(\varepsilon)}{t^{\lambda\gamma}}
\int\limits_{B(x,t)}|I_{\alpha}f(y)|^{q-\varepsilon}d\mu(y)\right)^{\frac{1}{q-\varepsilon}}
=\psi(\varepsilon)^{\frac{1}{q-\varepsilon}}\left(\frac{1}{t^{\lambda\gamma}}
\int\limits_{B(x,t)}|I_{\alpha}f(y)|^{q-\varepsilon}d\mu(y)\right)^{\frac{1}{q-\varepsilon}}
$$

$$
\leq(\psi(\varepsilon))^{\frac{1}{q-\varepsilon}}\sup_{\substack{x\in X\\
0\leq t<d}}\left(\frac{1}{t^{\lambda\gamma}}
\int\limits_{B(x,t)}|I_{\alpha}f(y)|^{q-\varepsilon}d\mu(y)\right)^{\frac{1}{q-\varepsilon}}
$$
$$
=\psi(\varepsilon)^{\frac{1}{q-\varepsilon}}\|I_{\alpha}f\|_{{\mathcal{L}}^{q-\varepsilon,\lambda}(X,\mu)}\leq $$

$$ \leq c(p-\eta,\alpha,\lambda)\psi(\varepsilon)^{\frac{1}{q-\varepsilon}} \|f\|_{{\mathcal{L}}^{q-\eta,\lambda}(X,\mu)}, $$

$$ \left(\text{where}\;\frac{1}{p-\eta}-\frac{1}{q-\varepsilon}=\frac{\alpha}{(1-\lambda)\gamma}\right) $$

$$ =c\frac{(1-\lambda)\gamma}{\alpha[(1-\lambda)\gamma-\alpha(p-\eta)]}\left(\left[\frac{p-\eta}{p-\eta-1}\right]
^{\frac{1}{q-\varepsilon}}+1\right)\psi(\varepsilon)^{\frac{1}{q-\varepsilon}}\|f\|_{{\mathcal{L}}^{p-\eta,\lambda}(X,\mu)}
$$

$$
\leq
\psi(\varepsilon)^{\frac{1}{q-\varepsilon}}\eta^{-\frac{\theta_1}{p-\eta}}
c(p-\eta,\alpha,\lambda)\eta^{\frac{\theta_1}{p-\eta}}\|f\|_{{\mathcal{L}}^{p-\eta,\lambda}(X,\mu)}
$$

$\Big($ observe that when $\sigma$ is small, then $\eta$ is also
small positive number; recall also that

$$ \varphi(u)=\left[p+\frac{(1-\lambda)(u-q)\gamma}{(1-\lambda)\gamma-\alpha(u-q)}\right]^
{\frac{\gamma(1-\lambda)-(u-q)\alpha}{(1-\lambda)\gamma}} \sim
u^{1+ \frac{\alpha q}{(1-\lambda)\gamma}}, \;\;\; u\to 0+, $$ and
$\psi(\varepsilon)^{\frac{1}{q-\varepsilon}}
\eta^{-\frac{\theta_1}{p-\eta}}=1$ $\Big)$

$$
\leq
\left[\sup\limits_{0<\eta\leq\sigma_{1}}c(p-\eta,\alpha,\lambda)\right]\|f\|_{{\mathcal{L}}^{p),\theta_1,
\lambda}(X,\mu)}
$$

Hence,
$$
\|I_{\alpha}f\|_{L^{q),\theta_2,\lambda}(X,\mu)}\leq\left[\sup\limits_{0<\eta\leq\sigma_{1}}
c(p-\eta,\alpha,\lambda)\right]\|f\|_{{\mathcal{L}}^{p),\theta_1,
\lambda}(X,\mu)},
$$
where
$$
c(p-\eta,\alpha,\lambda)=\overline{c}\frac{(1-\lambda)\gamma}{\alpha[(1-\lambda)\gamma-\alpha(p-\eta)]}
\bigg(
\left[\frac{p-\eta}{p-\eta-1}\right]^{\frac{1}{q-\varepsilon}}
+1\bigg),
$$
$\overline{c}$ is the constant independent of $p,\eta$ and
$\alpha$; $\sigma_{1}$ is a small positive number. Observe that if
$\sigma_{1}$ is sufficiently small, then
$(1-\lambda)\gamma-\alpha(p-\eta)\geq\eta_{0}>0$,
$\frac{p-\eta}{p-\eta-1} \leq p'+1$ for some $\eta_{0}$  when
$0<\eta\leq\sigma_{1}$.
\end{proof}

\subsection{Potential $
(T_{\alpha}f)(x)=\int\limits_{X}\frac{f(y)}{\mu
B(x,\rho(x,y))^{1-\alpha}}d\mu(y)$}
\setcounter{equation}{0}

Let

$$ (T_{\alpha}f)(x)=\int\limits_{X}\frac{f(y)}{\mu B(x,\rho(x,y))^{1-\alpha}}d\mu(y),\;\;0<\alpha<1. $$
Suppose that instead of condition (4.1) of the previous section  the following conditions hold:

(i)

$$ \mu \{ x\} =0, \;\;\;\; \text{for all }\; x\in X; $$

(ii)
$$
\mu B(x,t)\;\;\text{is continuous in}\;t\;\text{for every}\;x\in
X. \eqno{(4.4)}
$$
For example, (4.4) holds if $\mu\{y\in X:\;\;\rho(x,y)=t\}=0$ for
arbitrary $x\in X$ and $t\in [0,d)$.

We say that $f\in L^{p),\varphi(\cdot),\lambda}(X,\mu)$ if
$$
\|f\|_{L^{p),\varphi(\cdot),\lambda}(X,\mu)}:=\sup\limits_{0<\varepsilon<p-1}\sup_{\substack{x\in X\\
0\leq r<d}}\left(\frac{\varphi(\varepsilon)}{\mu
B(x,r)^{\lambda}}\int\limits_{B(x,r)}|f(y)|^{p}d\mu(y)\right)^{1/p}<\infty,\;\;
0\leq \lambda<1,
$$
where $\varphi$ is a positive function in $(0, p-1)$ which is
increasing near $0$ and satisfies the condition $\varphi (0+)=0$.
If $\varphi(\varepsilon)= \varepsilon^{\theta}$, where $\theta$ is
a positive number, then we denote $L^{p), \varphi(\cdot),
\lambda}(X, \mu)$ by $L^{p), \theta, \lambda}(X, \mu)$.

Our aim in this section is to prove the next statement:

%Theorem 3.4.1
{\bf Theorem 5.1.} {\em Let $1<p<\infty$,
$0<\alpha<\frac{1-\lambda}{p}$, $0\leq\lambda<1$ and let
$\frac{1}{p}-\frac{1}{q}=\frac{\alpha}{1-\lambda}$. Let $\theta_1$
be a positive number. We set
$$\theta_2= \theta_1\Big(1 + \frac{\alpha q}{1-\lambda}\Big).$$
Then the operator $T_{\alpha}$ is bounded from $L^{p),\theta_1,
\lambda}(X,\mu)$ to $L^{q),\theta_2, \lambda}(X,\mu)$.}

%where $$ \varphi(u)=\left[p+\frac{(1-\lambda)(u-q)}{1-\lambda-\alpha(u-q)}\right]^
%{\frac{1-\lambda-(u-q)\alpha}{1-\lambda}}. $$

To prove Theorem 5.1 we need the following lemma

\vskip+0.1cm

%Lemma 3.4.1.
{\bf Lemma 5.1.} {\em  Let $1<p<\infty$,
$0<\alpha<\frac{1-\lambda}{p}$,
$\frac{1}{p}-\frac{1}{q}=\frac{\alpha}{1-\lambda}$, where
$0\leq\lambda<1$. Then the inequality
$$
\|T_{\alpha}f\|_{L^{q,\lambda}(X,\mu)}\leq c(p,\alpha,\lambda)
\|f\|_{L^{p,\lambda}(X,\mu)}
$$
holds, where
$$ c(p,\alpha,\lambda)=c\left(C_{\alpha}+ \frac{p}{1-\lambda-\alpha
p}\right) \left[(p')^{1/q}+1\right] $$ and the positive constant
$c$ does not depend on $p$ and $\alpha$.}

\begin{proof}
Following the idea of Hedberg \cite{Hed} and taking into account the proof of $(4.2)$  we have that

$$ |T_{\alpha}f(x)|\leq c(p,\lambda,\alpha)(Mf)^{1-\frac{p\alpha}{1-\lambda}}(x)\|f\|_{L^{p,\lambda}(X,\mu)}^{\frac{\alpha
p}{1-\lambda}}, \eqno{(4.5)} $$
where
$$
c(p,\lambda,\alpha)=C_{\alpha}+\frac{p}{1-\lambda-\alpha p}
$$
and
$$
(Mf)(x)=\sup_{\substack{x\in X\\
0\leq r<d}}\frac{1}{\mu B(x,r)}\int\limits_{B(x,r)}|f(y)|d\mu(y).
$$
We set
$$
\overline{f}(x,t):=\frac{1}{\mu
B\left(x,\frac{\rho(x,t)}{\overline{\alpha}}\right)}\int\limits_{B\left(x,\frac{\rho(x,t)}{\overline{\alpha}}
\right)}|f(y)| d\mu(y), $$ where $\overline{\alpha}$ is the
constant between $0$ and $1$. In fact $\overline{\alpha}$ is the
constant from the reverse doubling condition $(3.3)$ (we use the
symbol $\overline{\alpha}$ instead of $\alpha$).

Observe that
$$ |I_{\alpha}f(x)|\leq
b\int\limits_{X}|f(y)|\left[\int\limits_{\overline{\alpha}\rho(x,y)<\rho(x,t)<\rho(x,y)}\mu
B(x,\rho(x,t))^{\alpha-2}d\mu(t)\right]d\mu(y),
$$
where $b$ is the constant depending on $\beta$ from $(3.3)$.

Hence,
$$ |I_{\alpha}f(x)|\leq b\int\limits_{X}\mu B(x,\rho(x,t))^{\alpha-2}
\left(\int\limits_{B\left(x,\frac{\rho(x,t)}{\overline{\alpha}}\right)}|f(y)|d\mu(y)\right)d\mu(t)
$$

$$ \leq b_{0}\int\limits_{X}\mu
B(x,\rho(x,t))^{\alpha-1}\overline{f}(t,x)d\mu(t), $$ where
$b_{0}$ is the positive constant which does not depend on $p$,
$\alpha$ and $\lambda$.

We take $\varepsilon>0$ which will be chosen later. Then
$$
|I_{\alpha}f(x)|\leq b_{0}\bigg[\int\limits_{B(x,\varepsilon)}\mu
B(x,\rho(x,t))^{\alpha-1}\overline{f}(t,x)d\mu(t)
$$
$$
+\int\limits_{X\backslash B(x,\varepsilon)}\mu
B(x,\rho(x,t))^{\alpha-1}\overline{f}(t,x)d\mu(t)\bigg]
=:b_{0}\left[J^{(1)}(x,t)+J^{(2)}(x,t)\right].
$$
It is easy  to see that (see also \cite{EdKoMe}, p. 348)
$$
J^{(1)}(x,t) \leq Mf(x) \int\limits_{B(x, \varepsilon)} \mu B(x,
\rho(x,t))^{\alpha-1} d\mu(t) \leq c_{\alpha}Mf(x)\mu
B(x,\varepsilon)^{\alpha},
$$
where the positive constant $c_{\alpha}$ depends only on $\alpha$.

Further, by H\"{o}lder's inequality we find that
$$ f(t,x)\leq\frac{1}{\mu B\left(x,\frac{\rho(x,t)}{\overline{\alpha}}\right)}
\left(\int\limits_{B\left(x,\frac{\rho(x,t)}{\overline{\alpha}}\right)}|f(t)|^{p}d\mu(t)\right)^{\frac{1}{p}}\mu
B\left(x,\frac{\rho(x,t)}{\overline{\alpha}}\right)^{\frac{1}{p'}}
$$
Besides this, by the inequality
$$ \int\limits_{X\backslash
B(x,\varepsilon)}\mu B(x,\rho(x,t))^{\alpha-1+\frac{\lambda-1}{p}}d\mu(t)\leq c \mu
B(x,\varepsilon)^{\alpha+\frac{\lambda-1}{p}}$$
(see Proposition 6.1.2 of \cite{EdKoMe}) we obtain that
$$
|I_{\alpha}f(x)|\leq b_{0}\bigg[ c_{\alpha}Mf(x)\mu
B(x,\varepsilon)^{\alpha}+\left(\int\limits_{X\backslash
B(x,\varepsilon)}\mu
B(x,\rho(x,t))^{\alpha-1+\frac{\lambda-1}{p}}d\mu(t)\right)
$$
$$
\|f\|_{L^{p,\lambda}(X,\mu)}\bigg]=b_{0}\bigg[c_{\alpha}(Mf)(x)\mu
B(x,\varepsilon)^{\alpha}+\overline{c}_{p,\lambda,\alpha}\mu
B(x,\varepsilon)^{\alpha+\frac{\lambda-1}{p}}
\|f\|_{L^{p,\lambda}(X,\mu)}\bigg], $$
where the positive constant $c_{\alpha}$ depends only on $\alpha$
and the positive constant $\overline{c}_{p,\lambda,\alpha}$ is
given by $\overline{c}_{p,\lambda,\alpha}=\frac{p}{1-\lambda-
\alpha p}$; $b_{0}$ does not depend on $p$, $\alpha$ and
$\lambda$.

Now we take (recall that $\mu B(x,\varepsilon)$ is continuous in
$\varepsilon$)
$$ \mu B(x,\varepsilon)=\left[\frac{\|f\|_{L^{p,\lambda}(X,\mu)}}{Mf(x)}\right]^{\frac{p}{1-\lambda}}. $$
Consequently,
$$
|I_{\alpha}f(x)|\leq\left(c_{\alpha}-\frac{p}{\alpha
p-1+\lambda}\right) \|f\|_{L^{p,\lambda}(X,\mu)}^{\frac{\alpha
p}{1-\lambda}}(Mf)^{1-\frac{p\alpha}{1-\lambda}}(x).
$$
Using the condition
$\frac{1}{p}-\frac{1}{q}=\frac{\alpha}{1-\lambda}$ and Proposition 2.2 we find that

$$ \left(\frac{1}{\mu
B(x,r)^{\lambda}}\int\limits_{B(x,r)}|I_{\alpha}f(x)|^{q}d\mu(x)\right)^{1/q}
$$

$$
\leq\mu B(x,r)^{-\lambda/q}\left(c_{\alpha}-\frac{p}{\alpha p-1
+\lambda}\right)\left[\int\limits_{B(x,r)}(Mf(y))^
{q\left[1-\frac{p\alpha}{1-\lambda}\right]}d\mu(y)\right]^{\frac{1}{q}}
\|f\|_{L^{p,\lambda}(X,\mu)}^{\frac{\alpha
p}{1-\lambda}}$$

$$ =\left(c_{\alpha}-\frac{p}{\alpha p- 1+\lambda}\right)\left[\frac{1}{\mu
B(x,r)^{\lambda/q}}\int\limits_{B(x,r)}(Mf(y))^{p}d\mu(y)\right]^{\frac{1}{q}}
\|f\|_{L^{p,\lambda}(X,\mu)}^{\frac{\alpha p}{1-\lambda}}$$

$$ \leq \left(c_{\alpha}-\frac{p}{\alpha p-
1+\lambda}\right)\|Mf\|_{L^{p,\lambda}(X,\mu)}^{p/q}\|f\|
_{L^{p,\lambda}(X,\mu)}^{\frac{\alpha p}{1-\lambda}} $$

%$$ \leq \text{ (by the inequality} \;\;$\|Mf\|_{L^{p,\lambda}(X,\mu)} \leq
%c_{0}b^{\lambda/p}\left[(p')^{1/p}+1\right]\|f\|_{L^{p, \lambda}(X, \mu)$)} $$

$$ \leq\left(c_{\alpha}-\frac{p}{\alpha p- 1+\lambda}\right) c_{0}b^{\lambda}
\left[ (p')^{1/q} + 1 \right] \|f\|_{L^{p,\lambda}(X,\mu)}. $$
\end{proof}
\vskip+0.3cm

{\em Proof of Theorem} 5.1. By using Lemma 5.1
and repeating the arguments of the proof of Theorem 4.1  we
conclude that Theorem 5.1 holds. Details are omitted. $\;\;\;
\Box$

\vskip+1cm

 {\bf Acknowledgement.} The author expresses his gratitude to
 Professor V. Kokilashvili for drawing his attention to the
 problems studied in this work and helpful comments and remarks.

%The author was partially
%supported by the Georgian National Science Foundation Grant
%(project numbers: No. GNSF/ST09/23/3-100 and No. GNSF/ST07/3-169).
\vskip+1cm

\vskip+0.5cm

Author's Address:

A. Meskhi: \

A. Razmadze Mathematical Institute, M. Aleksidze St.,  Tbilisi
0193, Georgia

Second Address: Department of Mathematics,  Faculty of Informatics
and Control Systems, Georgian Technical University, 77, Kostava
St., Tbilisi, Georgia.

e-mail:  meskhi@@rmi.acnet.ge \vskip+0.5cm

\end{document}